\documentclass{llncs}[10pt]
\usepackage{llncsdoc}
\usepackage{amsmath,graphicx}
\usepackage{amsfonts} 
\usepackage{amssymb}
\usepackage[latin1]{inputenc}
\title{Warped metrics for location-scale models}

\author{Salem Said, Yannick Berthoumieu}

\institute{Laboratoire IMS (CNRS - UMR 5218), Universit\'e de Bordeaux \\ \email{{salem.said;yannick.berthoumieu}@ims-bordeaux.fr}}

\usepackage[update,prepend]{epstopdf}
\RequirePackage{fix-cm}
\DeclareMathSizes{10}{10}{5}{5}

\begin{document}

\maketitle

\begin{abstract}
This paper argues that a class of Riemannian metrics, called warped metrics, plays a fundamental role in statistical problems involving location-scale models. The paper reports three new results\,: i) the Rao-Fisher metric of any location-scale model is a warped metric, provided that this model satisfies a natural invariance condition, ii) the analytic expression of the sectional curvature of this metric, iii) the exact analytic solution of the geodesic equation of this metric. The paper applies these new results to several examples of interest, where it shows that warped metrics turn location-scale models into complete Riemannian manifolds of negative sectional curvature. This is a very suitable situation for developing algorithms which solve problems of classification and on-line estimation. Thus, by revealing the connection between warped metrics and location-scale models, the present paper paves the way to the introduction of new efficient statistical algorithms. 
\end{abstract}

\section{Introduction\,: definition and two examples}
This paper argues that a class of Riemannian metrics, called warped metrics, is natural and useful to statistical problems involving location-scale models. A warped metric is defined as follows~\cite{petersen}. Let $M$ be a Riemannian manifold with Riemannian metric $ds^2_{M\,}$. Consider the manifold $\mathcal{M} = M \times (0\,,\infty)\,$, equipped with the Riemannian metric,
\begin{equation} \label{eq:warp}
  ds^2(z) = I_0(\sigma)\,d\sigma^2 \,+\, I_1(\sigma)\,ds^2_M(\bar{x})
\end{equation}
where each $z \in \mathcal{M}$ is a couple $(\bar{x}\,,\sigma)$ with $\bar{x} \in M$ and $\sigma \in (0\,,\infty)$. The Riemannian metric (\ref{eq:warp}) is called a warped metric on $\mathcal{M}\,$. The functions $I_0$ and $I_1$ have strictly positive values and are part of the definition of this metric.  

The main claim of this paper is that warped metrics arise naturally as Rao-Fisher metrics for a variety of location-scale models. Here, to begin, two examples of this claim are given. Example 1 is classic, while Example 2, to our knowledge, is new in the literature. As of now, the reader is advised to think of $\mathcal{M}$ as a statistical manifold, where $\bar{x}$ is a location parameter and $\sigma$ is either a scale parameter or a concentration parameter. \\[0.1cm]
\textbf{Example 1 (univariate normal model)\,:} let $M = \mathbb{R}$, with $ds^2_M(\bar{x}) = d\bar{x}^2$ the canonical metric of $\mathbb{R\,}$. If each $z = (\bar{x}\,,\sigma)$ in $\mathcal{M}$ is identified with the univariate normal density of mean $\bar{x}$ and standard deviation $\sigma\,$, then the resulting Rao-Fisher metric on $\mathcal{M}$ is given by~\cite{atkinson}
\begin{equation} \label{eq:example1}
  ds^2(z) = \sigma^{-2}\,d\sigma^2 \,+\, \frac{1}{2}\sigma^{-2}\, d\bar{x}^2
\end{equation}
\textbf{Example 2 (von Mises-Fisher model)\,:} let $M = S^{\,2\,}$, the unit sphere with $ds^2_M = d\theta^{2}$ its canonical metric induced from $\mathbb{R}^{3\,}$. Identify $z = (\bar{x}\,,\sigma)$ in $\mathcal{M}$ with the von Mises-Fisher density of mean direction $\bar{x}$ and concentration parameter $\sigma\,$~\cite{mardia}. The resulting Rao-Fisher metric on $\mathcal{M}$ is given by
\begin{equation} \label{eq:example2}
  ds^2(z) = \left(\, \sigma^{-2} - \sinh^{-2}\,\sigma\,\right)\,d\sigma^2 \,+\, \left(\, \sigma\,\coth\,\sigma - 1 \,\right)\, d\theta^2(\bar{x})
\end{equation}
\textbf{Remark a\,:} note that $\sigma$ is a scale parameter in Example 1, but a concentration parameter in Example 2. Accordingly, at $\sigma = 0$, the metric (\ref{eq:example1}) becomes infinite, while the metric (\ref{eq:example2}) remains finite and degenerates to $\left.ds^2(z)\,\right|_{\sigma = 0}\, = (1/3)\,d\sigma^2$. Thus, (\ref{eq:example2}) gives a Riemannian metric on the larger Riemannian manifold $\hat{\mathcal{M}} = \mathbb{R}^3$, which contains $\mathcal{M}$, obtained by considering $\sigma$ as a radial coordinate and $\sigma = 0$ as the origin of $\mathbb{R}^3$. \hfill $\blacksquare$

\section{A general theorem : from Rao-Fisher to warped metrics}
Examples 1 and 2 of the previous section are special cases of Theorem 1, given here. To state this theorem, let $(M,ds^2_M)$ be an irreducible Riemannian homogeneous space, under the action of a group of isometries $G$~\cite{kobnom}. Denote by $g\cdot x$ the action of $g \in G$ on $x \in M$. Then, assume each $z = (\bar{x}\,,\sigma)$ in $\mathcal{M}$ can be identified uniquely and regularly with a probability density $p(x|\,z)=p(x|\,\bar{x}\,,\sigma)$ on $M$, with respect to the Riemannian volume element, such that the following property is verified,
\begin{equation} \label{eq:invarianceh}
  p(\;\!g\cdot x|\,g\cdot \bar{x}\,,\sigma) \,= \, p(x|\,\bar{x}\,,\sigma) \hspace{0.3cm} g \in G
\end{equation}
The densities $p(x|\,\bar{x}\,,\sigma)$ form a statistical model on $M$, where $\bar{x}$ is a location parameter and $\sigma$ can be chosen as either a scale or a concentration parameter, (roughly, a scale parameter is the inverse of a concentration parameter). 

In the statement of Theorem 1, $\ell(z) = \log\,p(x|\,z)$ and $\nabla_{\bar{x}}\,\ell(z)$ denotes the Riemannian gradient vector field of $\ell(z)$, with respect to $\bar{x} \in M$. Moreover, $\left\Vert\nabla_{\bar{x}}\ell(z)\right\Vert$ denotes the length of 
this vector field, as measured by the metric $ds^2_{M\,}$. 
\begin{theorem}[warped metrics]
  The Rao-Fisher metric of the statistical model $\left\lbrace\, p(x|\,z) \, ; \, z \in \mathcal{M}\,\right\rbrace$ is a warped metric of the form (\ref{eq:warp}), defined by
\begin{equation} \label{eq:I}
  I_0(\sigma) = \mathbb{E}_{z} \left(\, \partial_\sigma\ell(z)\,\right)^2 \hspace{0.3cm}
  I_1(\sigma) = \left.\mathbb{E}_{z}\left\Vert\nabla_{\bar{x}}\ell(z)\right\Vert^2\middle/\mathrm{dim}\,M\right.
\end{equation}
where  $\mathbb{E}_{z}$ denotes expectation with respect to $p(x|\,z)\,$. Due to property (\ref{eq:invarianceh}), the two expectations appearing in (\ref{eq:I}) do not depend on the parameter $\bar{x}$, so $I_0$ and $I_1$ are well-defined functions of $\sigma$.  
\end{theorem}
\textbf{Remark b\,:} the proof of Theorem 1 cannot be given here, due to lack of space. It relies strongly on the assumption that the Riemannian homogeneous space $M$ is irreducible. In particular, this allows the application of Schur's lemma, from the theory of group representations~\cite{chevalley}. To say that $M$ is an irreducible Riemannian homogeneous space means that the following property is verified\,: if $K_{\bar{x}}$ is the stabiliser in $G$ of $\bar{x} \in M\,$, then the isotropy representation $k \mapsto \left. dk\right|_{\bar{x}}$ is an irreducible representation of $K_{\bar{x}}$ 
in the tangent space $T_{\bar{x}}M\,$. \hfill $\blacksquare$ \\[0.1cm]
\textbf{Remark c\,:} if the assumption that $M$ is irreducible is relaxed, then Theorem 1 generalises to a similar statement, involving so-called multiply warped metrics. Roughly, this is because a homogeneous space which is not irreducible, may still decompose into a direct product of irreducible homogeneous spaces~\cite{kobnom}. \hfill $\blacksquare$ \\[0.1cm] 
\textbf{Remark d\,:} statistical models on $M$ which verify (\ref{eq:invarianceh}) often arise under an exponential form,
\begin{equation} \label{eq:exponential}
 p(x|\,\bar{x}\,,\sigma) \,= \, \exp\left(\, \eta\cdot D(x\,,\bar{x}) - \psi(\eta)\,\right)
\end{equation}
where $\eta = \eta(\sigma)$ is a natural parameter, and $\psi(\eta)$ is the cumulant generating function of the statistic $D(x\,,\bar{x})\,$. Then, for assumption (\ref{eq:invarianceh}) to hold, it is necessary and sufficient that
\begin{equation} \label{eq:invarianceD}
   D(\;\!g\cdot x\,,\,g\cdot\bar{x})\,=\, D(x\,,\bar{x})
\end{equation}
Both examples 1 and 2 are of the form (\ref{eq:exponential}), as is Example 3, in the following section, which deals with the Riemannian Gaussian model~\cite{said1}\cite{said2}.  \hfill $\blacksquare$\\[0.1cm]

\section{Curvature equations and the extrinsic geometry of $M$}
For each $\sigma \in (0\,,\infty)\,$, there is an embedding of $M$ into $\mathcal{M}\,$, as the surface $M \times \lbrace \sigma \rbrace\,$. This embedding yields an extrinsic geometry of $M$, given by the first and second fundamental forms~\cite{docarmo}. 

The first fundamental form is the restriction of the metric $ds^2$ of $\mathcal{M}\,$ to the tangent bundle of $M$. This will be denoted $ds^2_M(x|\,\sigma)$ for $x \in M$. It is clear from (\ref{eq:warp}) that
\begin{equation} \label{eq:firstform}
  ds^2_M(x|\,\sigma) \,=\, I_1(\sigma)\, ds^2_M(x)
\end{equation}
This extrinsic Riemannian metric on $M$ is a scaled version of its intrinsic metric $ds^2_{M\,}$. It induces an extrinsic Riemannian distance given by
\begin{equation} \label{eq:extrinsicdistance}
  d^{\,2}(x\,,y|\,\sigma) \,=\, I_1(\sigma)\, d^{\,2}(x\,,y) \hspace{0.3cm} x\,,\, y \, \in M
\end{equation}
where $d(x\,,y)$ is the intrinsic Riemannian distance, induced by the metric $ds^2_M\,$. 

The extrinsic distance (\ref{eq:extrinsicdistance}) is a generalisation of the famous Mahalanobis distance. In fact, replacing in Example 1 yields the classical expression of the Mahalanobis distance $d^{\,2}(x\,,y) \,= |x-y|^2/2\sigma^2$. The significance of this distance can be visualised as follows\,: if $\sigma$ is a dispersion parameter, the extrinsic distance between two otherwise fixed points $x\,,y \,\in M$ will decrease as $\sigma$ increases, as if the space $M$ were contracting, (for a concentration parameter, there is an expansion, rather than a contraction).

The second fundamental form is given by the tangent component of the covariant derivative of the unit normal to the surface $M \times \lbrace \sigma \rbrace\,$. This unit normal is $\partial_r$ where $r$ is the vertical distance coordinate, given by $dr/d\sigma = I^{\frac{1}{2}}_0(\sigma)\,$. Using Koszul's formula~\cite{helgason}, it is possible to express the second fundamental form,
\begin{equation} \label{eq:secondform}
  S(v) \,=\, \frac{1}{2}\left(\,\partial_r I_1\middle/ I_1\,\right)\,v 
\end{equation}
for any $v$ tangent to $M$. Knowledge of the second fundamental form is valuable, as it yields the relationship between extrinsic and intrinsic curvatures of $M$. 
\begin{proposition}[curvature equations]
Let $K^\mathcal{M}$ and $K^M$ denote the sectional curvatures of $\mathcal{M}$ and $M$\,. The following are true
\begin{equation} \label{eq:gauss}
   K^\mathcal{M}(u\,,v) \,=\, \left(\,1\middle/ I_1\,\right)\,K^M(u\,,v) \,-\, \frac{1}{4}\left(\,\partial_r I_1\middle/ I_1\,\right)^2
\end{equation}
\begin{equation} \label{eq:mixed}
   K^\mathcal{M}(u\,,\partial_r\,) \,=\, -\, \left(\, \partial^{\,2}_r \, I^{\frac{1}{2}}_1\middle/I^{\frac{1}{2}}_1 \,\right)
\end{equation}
for any linearly independent $u\,,v$ tangent to $M\,$. 
\end{proposition}
\textbf{Remark e\,:} here, Equation (\ref{eq:gauss}) is the Gauss curvature equation. Roughly, it shows that embedding $M$ into $\mathcal{M}$ adds negative curvature. Equation (\ref{eq:mixed}) is the mixed curvature equation. If the intrinsic sectional curvature $K^M$ is negative, then (\ref{eq:gauss}) and (\ref{eq:mixed}) show that the sectional curvature $K^\mathcal{M}$ of $\mathcal{M}$ is negative if and only if $I^{\frac{1}{2}}_1$ is a convex function of the vertical distance $r$\,.  \hfill$\blacksquare$ \\[0.1cm]
 \textbf{Return to example 1\,:} here, $M = \mathbb{R}$ is one-dimensional, so the Gauss equation (\ref{eq:gauss}) does not provide any information. The mixed curvature equation gives the curvature of the two-dimensional manifold $\mathcal{M}$. In this equation, $\partial_r = \sigma\, \partial_\sigma\,$, and it follows that
\begin{equation} \label{eq:curvature1}
 K^\mathcal{M}(u\,,\partial_r\,) \,=\, -\, 1
\end{equation}
so $\mathcal{M}$ has constant negative curvature. In fact, it was observed long ago that the metric (\ref{eq:example1}) is essentially the Poincar\'e half-plane metric~\cite{atkinson}. \hfill $\blacksquare$ \\[0.1cm]
 \textbf{Return to example 2\,:} in Example 2, $M = S^2$  so $K^M \equiv 1$ is constant. It follows from the Gauss equation that each sphere $S^2 \times \lbrace \sigma \rbrace$ has constant extrinsic curvature, equal to
\begin{equation} \label{eq:curvature2}
  \left.K^\mathcal{M}\right|_\sigma \,=\, \, \left(\,1\middle/ I_1\,\right)\, \,-\, \frac{1}{4}\left(\,\partial_r I_1\middle/ I_1\,\right)^2
\end{equation}
Upon replacing the expressions of $I_1$ and $\partial_r$ based on (\ref{eq:example2}), this is found to be strictly negative for $\sigma > 0$,
\begin{equation} \label{eq:crazy}
\left.K^\mathcal{M}\right|_\sigma \,< 0 \hspace{0.3cm} \text{ for } \;\sigma > 0
\end{equation}
 Thus, the Rao-Fisher metric (\ref{eq:example2}) induces a negative extrinsic curvature on each spherical surface $S^2 \times \lbrace \sigma \rbrace$\,. In fact, by studying the mixed curvature equation (\ref{eq:mixed}), it is seen the whole manifold $\mathcal{M}$ equipped with the Rao-Fisher metric (\ref{eq:example2}) is a manifold of negative sectional curvature.   \hfill$\blacksquare$ \\[0.1cm]
\textbf{Example 3 (Riemannian Gaussian model)\,:} a Riemannian Gaussian distribution may be defined on any Riemannian symmetric space $M$ of non-positive curvature. It is given by the probability density with respect to Riemannian volume
\begin{equation} \label{eq:rgd}
  p(x|\,\bar{x}\,,\sigma) = Z^{-1}(\sigma)\,\exp\left[\, - \frac{d^{\,2}(x\,,\bar{x})}{2\sigma^2}\,\right]
\end{equation}
where the normalising constant $Z(\sigma)$ admits a general expression, which was given in~\cite{said2}. If $M$ is an irreducible Riemannian symmetric space, then Theorem 1 above applies to the Riemannian Gaussian model (\ref{eq:rgd}), leading to a warped metric with
\begin{equation} \label{eq:example3}
   I_0(\sigma) \,= \, \psi^{\prime\prime}(\eta) \hspace{0.3cm} I_1(\sigma) \,= \, \left. 4\eta^2\, \psi^{\prime}(\eta)\middle/\mathrm{dim}\,M\right.
\end{equation}
where $\eta = -1/2\sigma^2$ and $\psi(\eta) = \log\, Z(\sigma)$. The result of equation (\ref{eq:example3}) is here published for the first time. Consider now the special case where $M$ is the hyperbolic plane. The analytic expression of $I_0$ and $I_1$ can be found from (\ref{eq:example3}) using
\begin{equation} \label{eq:Z}
  Z(\sigma) \,=\, \mathrm{Const}. \,\, \sigma\,\times\, e^{\,\sigma^2/4}\,\times\,\mathrm{erf}(\sigma/2)
\end{equation}
which was derived in~\cite{said1}. Here, $\mathrm{erf}$ denotes the error function. Then, replacing (\ref{eq:example3}) in the curvature equations (\ref{eq:gauss}) and (\ref{eq:mixed}) yields the same result as for Example 2\,: the manifold $\mathcal{M}$ equipped with the Rao-Fisher metric (\ref{eq:example3}) is a manifold of negative sectional curvature. \hfill$\blacksquare$ \\[0.1cm]
\textbf{Remark f (a conjecture)\,:} based on the three examples just considered,  it seems reasonable to conjecture that warped metrics arising from Theorem 1 will always lead to manifolds $\mathcal{M}$ of negative sectional curvature. \hfill$\blacksquare$ 

\section{Solution of the geodesic equation\,: conservation laws}
If the assumptions of Theorem 1 are slightly strengthened, then an analytic solution of the geodesic equation of the Riemannian metric (\ref{eq:warp}) on $\mathcal{M}$ can be obtained, by virtue of the existence of a sufficient number of conservation laws. To state this precisely, let $\langle\cdot,\cdot\rangle_M$ and $\langle\cdot,\cdot\rangle_\mathcal{M}$ denote respectively the scalar products defined by the metrics $ds^2_M$ and $ds^2$. 

Two kinds of conservation laws hold along any affinely parameterised geodesic curve $\gamma(t)$ in $\mathcal{M}$, with respect to the metric $ds^2\,$. These are conservation of energy and conservation of moments~\cite{gallot}. If the geodesic $\gamma(t)$ is expressed as a couple $(\,\sigma(t)\,,x(t)\,)$ where $\sigma(t) > 0$ and $x(t) \in M$\,, then the energy of this geodesic is 
\begin{equation} \label{eq:energy}
   E\,=\, I_0(\sigma)\,\dot{\sigma}^2 \,+\, I_1(\sigma)\,\Vert \dot{x}\Vert^2
\end{equation}
where the dot denotes differentiation with respect to $t\,$, and $\Vert \dot{x}\Vert$ the Riemannian length of $\dot{x}$ as measured by the metric $ds^2_{M\,}$.

 On the other hand, if $\xi$ is any element of the Lie algebra of the group of isometries $G$ acting on $M$, the corresponding moment of the geodesic $\gamma(t)$ is
\begin{equation} \label{eq:moment}
   J(\xi) \,=\, I_1(\sigma)\,\langle\, \dot{x}\,,X_\xi\,\rangle_M
\end{equation}
where $X_\xi$ is the vector field on $M$ given by $X_\xi(x) = \left.\frac{d}{dt}\right|_{t=0}\,e^{t\xi}\cdot x\,$. The equation of the geodesic $\gamma(t)$ is given as follwos.
\begin{proposition}[conservation laws and geodesics]
  For any geodesic $\gamma(t)$, its energy $E$ and its moment $J(\xi)$ for any $\xi$ are conserved quantities, remaining constant along this geodesic. If $M$ is an irreducible Riemannian symmetric space, the equation of the geodesic $\gamma(t)$ is the following,
\begin{equation} \label{eq:xt}
   x(t)\,=\, \mathrm{Exp}_{\,x(0)}\left[\,\left(\,\int^t_0\, \frac{I_1(\sigma(0))}{I_1(\sigma(s))}ds \,\right)\dot{x}(0)\,\right]
\end{equation}
\begin{equation} \label{eq:sigmat}
  t \,=\,\pm\,\int^{\sigma(t)}_{\sigma(0)}\, \frac{I^{\frac{1}{2}}_0(\sigma)\,d\sigma}{\sqrt{E - V(\sigma)}} 
\end{equation}
where $\mathrm{Exp}$ denotes the Riemannian exponential mapping of the metric $ds^2_M$ on $M\,$, and $V(\sigma)$ is the function
$V(\sigma) = J_0\times \left.I_1(\sigma(0))\middle/I_1(\sigma)\right.\,$, with $J_0 = I_1(\sigma(0))\,\Vert\dot{x}(0)\Vert^{2\,}$.
\end{proposition}
\textbf{Remark g\,:} under the assumption that $M$ is an irreducible Riemannian symmetric space, the second part of Proposition 2, stating the equations of $x(t)$ and $\sigma(t)$ is a corollary of the first part, stating the conservation of energy and moment. The proof, as usual not given due to lack of space, relies on a technique of lifting the geodesic equation to the Lie algebra of the group of isometries $G$. \hfill$\blacksquare$\\[0.1cm]
\textbf{Remark h\,:} here, Equation (\ref{eq:xt}) states that $x(t)$ describes a geodesic curve in the space $M$, with respect to the metric $ds^2_{M\,}$, at a variable speed equal to $\left.I_1(\sigma(0))\middle/I_1(\sigma(t))\right.$. Equation (\ref{eq:sigmat}) states that $\sigma(t)$ describes the one-dimensional motion of a particle of energy $E$ and mass $2I_0(\sigma)$, in a potential field $V(\sigma)$. \hfill$\blacksquare$ \\[0.1cm]
\textbf{Remark i (completeness of $\mathcal{M}$)\,:} from Equation (\ref{eq:sigmat}) it is possible to see that any geodesic $\gamma(t)$ in $\mathcal{M}$ is defined for all $t > 0$, if and only if the following conditions are verified
\begin{equation} \label{eq:completeness}
   \int_0 \,  I^{\frac{1}{2}}_0(\sigma)\,d\sigma = \infty \hspace{0.6cm}
 \int^\infty \,  I^{\frac{1}{2}}_0(\sigma)\,d\sigma = \infty
\end{equation}
where the missing integration bounds are arbitrary. The first condition ensures that $\gamma(t)$ may not escape to $\sigma = 0$ within a finite time, while the second condition ensures the same for $\sigma = \infty$. The two conditions (\ref{eq:completeness}), taken together, are necessary and sufficient for $\mathcal{M}$ to be a complete Riemannian manifold. \hfill $\blacksquare$ \\[0.1cm]
\textbf{Return to Example 2\,:} for the von Mises-Fisher model of Example 2, the second condition in (\ref{eq:completeness}) is verified, but not the first. Therefore, a geodesic $\gamma(t)$ in $\mathcal{M}$ may escape to $\sigma = 0$ within a finite time. However, $\gamma(t)$ is also a geodesic in the larger manifold $\hat{\mathcal{M}} = \mathbb{R}^{3\,}$, which contains $\sigma = 0$ as its origin. If $\gamma(t)$ arrives at $\sigma = 0$ at some finite time, it will just go through this point and immediately return to $\mathcal{M}$. In fact, $\hat{\mathcal{M}}$ is a complete Riemannian manifold which has $\mathcal{M}$  as an isometrically embedded submanifold. \hfill $\blacksquare$
\section{The road to applications: classification and estimation}
The theoretical results of the previous chapters have established that warped metrics are natural statistical objects arising in connection with location-scale models, which are invariant under some group action. Precisely, Theorem 1 has stated that warped metrics appear as Rao-Fisher metrics for all location-scale models which verify the group invariance condition (\ref{eq:invarianceh}). 

Analytical knowledge of the Rao-Fisher metric of a statistical model is potentially useful to many applications. In particular, to problems of classification and efficient on-line estimation. However, in order for such applications to be realised, it is necessary for the Rao-Fisher metric to be well-behaved. Propositions 2 and 3 in the above seem to indicate such a good behavior for warped metrics on location-scale models. 

Indeed, as conjectured in Remark f, the curvature equations of Proposition 2 would indicate that the sectional curvature of these warped metrics is always negative. Then, if the conditions for completeness, given in Remark i based on Proposition 3, are verified, the location-scale models equipped with these warped metrics appear as complete Riemannian manifolds of negative curvature. This is a favourable scenario, (which at least holds for the von Mises-Fisher model of Example 2), under which many algorithms can be implemented. 

For classification problems, it becomes straightforward to find the analytic expression of Rao's Riemannian distance, and to compute Riemannian centres of mass, whose existence and uniqueness will be guaranteed. These form the building blocks of many classification methodologies. 

For efficient on-line estimation, Amari's natural gradient algorithm turns out to be identical to the stochastic Riemannian gradient algorithm, defined using the Rao-Fisher metric. Then, analytical knowledge of the Rao-Fisher metric, (which is here a warped metric), and of its completeness and curvature properties, yields an elegant formulation of the natural gradient algorithm, and a geometrical means of proving its efficiency and understanding its convergence properties.


\bibliographystyle{splncs}
\scriptsize
\bibliography{refs}

\end{document}